\documentclass[11pt]{article}

\usepackage{latexsym}
\usepackage{amssymb}

\newtheorem{theorem}{Theorem}
\newtheorem{proposition}{Proposition}
\newtheorem{corollary}{Corollary}

\newtheorem{lemma}{Lemma}

\newtheorem{example}{Example}

\begin{document}
{
\begin{center}
{\Large\bf
On some Sobolev spaces with matrix weights and classical type Sobolev orthogonal polynomials.}
\end{center}
\begin{center}
{\bf S.M. Zagorodnyuk}
\end{center}

\section{Introduction.}

The theory of Sobolev orthogonal polynomials attracted a lot of attention in the past 30 years, see surveys in~\cite{cit_5150_M_X_Survey_2015},
\cite{cit_5150_Survey_2006}, \cite{cit_5150__Survey_2001}, \cite{cit_5150__Survey_1998}, \cite{cit_5150__Survey_1996}, 
\cite{cit_5150_Survey_1993}. It is still under development and many aspects (such as the existence of recurrence relations) are 
hidden (\cite{cit_2000__Kwon___2009},\cite{cit_1500___Kim__2014}).
It turned out that it is convenient to construct new families of Sobolev orthogonal polynomials by using known orthogonal polynomials
on the real line (OPRL) or orthogonal polynomials on the unit circle (OPUC), 
see~\cite{cit_9000__Zagorodnyuk_JAT_2020},\cite{cit_9000__Zagorodnyuk_CMA_2020}.
Moreover, if a system of OPRL or OPUC is an eigenvector of a pencil of differential equations, then the associated Sobolev orthogonal
polynomials have a similar property as well. As for the existence of a recurrence relation, this question is
more complicated and needs additional efforts.
Let $\mathcal{K}$ denote the real line or the unit circle. The following problem seems to be an appropriate framework to
study classical type Sobolev orthogonal polynomials (cf.~\cite[Problem 1]{cit_9000__Zagorodnyuk_CMA_2020}).

\noindent
\textbf{Problem 1.} 
\textit{
To describe all Sobolev orthogonal polynomials $\{ y_n(z) \}_{n=0}^\infty$ on $\mathcal{K}$, satisfying the following two properties:
\begin{itemize}
 \item[(a)]  Polynomials $y_n(z)$ satisfy the following differential equation:
\begin{equation}
\label{f1_42}
R y_n(z) = \lambda_n S y_n(z),\qquad n=0,1,2,...,
\end{equation}
where $R,S$ are linear differential operators of finite orders, having complex polynomial coefficients not depending on $n$;
$\lambda_n\in\mathbb{C}$;
 \item[(b)]  Polynomials $y_n(z)$ obey the following difference equation:
\begin{equation}
\label{f1_43}
L \vec y(z) = z M \vec y(z),\quad \vec y(z) = (y_0(z),y_1(z),...)^T,
\end{equation}
where $L,M$ are semi-infinite complex banded (i.e. having a finite number of non-zero diagonals) matrices. 
\end{itemize}
}

Relation~(\ref{f1_42}) means that $y_n(z)$ are eigenvalues of the operator pencil $R-\lambda S$, and relation~(\ref{f1_43})
shows that vectors of $y_n(z)$ are eigenvalues of the operator pencil $L - z M$.
For a background on operator pencils (or operator polynomials) see~\cite{cit_7000_Markus},\cite{cit_8000_Rodman}.

In the case: $y_n(z) = z^n$, $\mathcal{K}=\mathbb{T}=\{ z\in\mathbb{C}:\ |z|=1 \}$, we have the following differential equation:
\begin{equation}
\label{f1_50}
z (z^n)' = n z^n,\qquad n\in\mathbb{Z}_+. 
\end{equation}
On the other hand, $y_n$ satisfies~(\ref{f1_43}) with $L$ being the identity semi-infinite matrix, $M$ being the semi-infinite matrix with all $1$
on the first subdiagonal and $0$ on other places.
We should emphasize that in Problem~1 we do not exclude OPRL or OPUC. They are formally considered as Sobolev orthogonal
polynomials with the derivatives of order $0$. In this way, we may view systems from Problem~1 as generalizations of systems
of classical orthogonal polynomials (see, e.g., \cite{cit_5105_Koekoek_book}).

Let us briefly describe the content of the paper.
In Section~2 we shall formulate the scheme from~\cite{cit_9000__Zagorodnyuk_JAT_2020} in a more general setting. 
We define a Sobolev space corresponding to a matrix measure and consider measurable factorizations of the corresponding matrix
weight. These factorizations can have several applications. The main application is to the construction of new Sobolev orthogonal 
polynomials. 
Some transparent conditions 
for the effectiveness of this scheme are given in
Proposition~\ref{p2_1}. 
Secondly, an application to the question of the density of polynomials in the
associated Sobolev space can be given (Corollary~\ref{c2_1}).
Finally, measurable factorizations can lead to new measures of orthogonality for Sobolev orthogonal polynomials which
are discussed at the end of Section~2.

In Section~3 we shall construct two concrete families of Sobolev orthogonal polynomials (depending on an arbitrary number of
complex parameters) which satisfy all conditions of Problem~1, see Corollaries~\ref{c3_3} and~\ref{c3_4}.
These polynomials possess explicit integral representations as well as explicit orthogonality relations (Theorem~\ref{t3_3}).
In particular cases, additional information on the location of zeros and asymptotics is given (Theorem~\ref{t3_1}).
For any system $\{ p_n \}_{n=0}^\infty$ of OPRL or OPUC we construct Sobolev orthogonal polynomials with explicit integral representations
involving $p_n$ (Corollaries~\ref{c3_1},\ref{c3_2}).

\noindent
{\bf Notations. }
As usual, we denote by $\mathbb{R}, \mathbb{C}, \mathbb{N}, \mathbb{Z}, \mathbb{Z}_+$,
the sets of real numbers, complex numbers, positive integers, integers and non-negative integers,
respectively. 
By $\mathbb{Z}_{k,l}$ we mean all integers $j$ satisfying the following inequality:
$k\leq j\leq l$; ($k,l\in\mathbb{Z}$).
For a complex number $c$ we denote
$[c]_0 = 1$, $[c]_k = c(c-1)...(c-k+1)$, $k\in\mathbb{N}$.
For $a\in\mathbb{R}$, we denote by $[a]$ the greatest integer number $n$: $n\leq a$.
By $\mathbb{C}_{m\times n}$ we denote the set of all $(m\times n)$ matrices with complex entries,
$\mathbb{C}_{n} := \mathbb{C}_{1\times n}$, $m,n\in\mathbb{N}$.
By $\mathbb{C}_{n\times n}^\geq$ we denote the set of all nonnegative Hermitian matrices from $\mathbb{C}_{n\times n}$,
$n\in\mathbb{N}$.
For $A\in\mathbb{C}_{m\times n}$ the notation $A^*$ stands for the adjoint matrix ($m,n\in\mathbb{N}$),
and $A^T$ means the transpose of $A$.
Set $\mathbb{T} := \{ z\in\mathbb{C}:\ |z|=1 \}$. 
By $\mathbb{P}$ we denote the set of all polynomials with complex coefficients.
For an arbitrary Borel subset $K$ of the complex plane we denote by $\mathfrak{B}(K)$ the set of all Borel subsets of $K$.
Let $\mu$ be an arbitrary (non-negative) measure on $\mathfrak{B}(K)$.
By $L^2_\mu = L^2_{\mu,K}$ we denote the usual space of (the classes of the equivalence of) complex Borel measurable functions $f$ on $K$ such that
$\| f \|_{L^2_{\mu,K}}^2 := \int_K |f(z)|^2 d\mu < \infty$.

By
$(\cdot,\cdot)_H$ and $\| \cdot \|_H$ we denote the scalar product and the norm in a Hilbert space $H$,
respectively. The indices may be omitted in obvious cases.
For a set $M$ in $H$, by $\overline{M}$ we mean the closure of $M$ in the norm $\| \cdot \|_H$.

\section{Measurable matrix factorizations and new measures of orthogonality for
Sobolev orthogonal polynomials.}

Fix an arbitrary Borel subset $K$ of the complex plane and an arbitrary $\rho\in\mathbb{N}$. 
Let $M(\delta) = ( m_{k,l}(\delta) )_{k,l=0}^\rho$ be a $\mathbb{C}_{(\rho+1)\times (\rho+1)}^\geq$-valued function on 
$\mathfrak{B}(K)$, which entries are countably additive on $\mathfrak{B}(K)$ ($\delta\in\mathfrak{B}(K)$). The function $M(\delta)$ is 
said to be \textit{a non-negative Hermitian-valued measure on $(K,\mathfrak{B}(K))$}, see~\cite[p. 291]{cit_7000__Rosenberg_1964}.
In what follows we shall need the space $L^2(M)$ of $\mathbb{C}_{\rho+1}$-valued functions which are square-integrable
with respect to $M$. Let us recall its definition from~\cite{cit_7000__Rosenberg_1964}.

Denote by $\tau(\delta)$ \textit{the trace measure}, $\tau(\delta) := \sum_{k=0}^\rho m_{k,k}(\delta)$, $\delta\in\mathfrak{B}(K)$.
By $M'_\tau := dM/d\tau = ( dm_{k,l}/ d\tau )_{k,l=0}^{\rho}$, we denote \textit{the trace derivative of $M$}.
One means by $L^2(M)$ a set of all (classes of the equivalence of)
measurable vector-valued functions
$\vec f(z): K\rightarrow \mathbb{C}_{\rho+1}$, 
$\vec f = (f_0(z),f_1(z),\ldots,f_\rho(z))$, such that
$$ \| \vec f \|^2_{L^2(M)} := \int_K  \vec f(z) M'_\tau(z) \vec f^*(z) d\tau  < \infty. $$
Two functions $\vec f$ and $\vec u$ belong to the same class of the equivalence if and only if $\| \vec f - \vec u \|_{L^2(M)} = 0$.
It is known that $L^2(M)$ is a Hilbert space with the following scalar product:
\begin{equation}
\label{f2_3}
( \vec f, \vec g )_{L^2(M)} := \int_K  \vec f(z) M'_\tau(z) \vec g^*(z) d\tau,\qquad  \vec f,\vec g\in L^2(M). 
\end{equation}
As usual, we shall often work with the representatives instead of the corresponding classes of the equivalence.
It is known (see~\cite[p. 294]{cit_7000__Rosenberg_1964}, \cite[Lemma 2.1]{cit_8000__Weron_1974})
that one can consider an arbitrary $\sigma$-finite (non-negative) measure $\mu$, with respect to which all $m_{k,l}$ are
absolutely continuous, and set $M_0(z) = M_{0,\mu}(z) := dM/d\mu$ (the Radon-Nikodym derivative of $M$ with respect to $\mu$).
Then the integral in~(\ref{f2_3}) exists if and only if the following integral exists:
\begin{equation}
\label{f2_4}
\int_K  \vec f(z) M_{0,\mu}(z) \vec g^*(z) d\mu. 
\end{equation}
If the integrals exist, they are equal. Such measures $\mu$ we shall call \textit{admissible}. The matrix function $M_{0,\mu}(z)$
is said to be \textit{the (matrix) weight, corresponding to an admissible measure $\mu$}.

Denote by $A^2(M)$ a linear manifold in $L^2(M)$ including those classes of the equivalence $[\cdot]$ which possess 
a representative of the following form: 
\begin{equation}
\label{f2_5}
\vec f(z) = (f(z),f'(z),...,f^{(\rho)}(z)). 
\end{equation}
By $W^2(M)$ we denote the closure of $A^2(M)$ in the norm of $L^2(M)$. 
The subspace $W^2(M)$ is said to be \textbf{the Sobolev space with the matrix measure $M$}.
Elements of $A^2(M)$ will be also
denoted by their first components. Thus, we can write $f(z)$ instead of $\vec f(z)$ for an element in~(\ref{f2_5}), if this cause no misunderstanding.

Let $\mu$ be an admissible measure and $M_0(z)$ be the corresponding matrix weight.
Suppose that $1,z,z^2,...$, all belong to $W^2(M)$. 
We shall also assume that
\begin{equation}
\label{f2_7}
(p,p)_{W^2(M)} > 0,
\end{equation}
for an arbitrary non-zero $p\in\mathbb{P}$.
Then one can apply the Gram-Schmidt orthogonalization process to construct a system
$\{ y_n(z) \}_{n=0}^\infty$, $\deg y_n=n$, of Sobolev orthogonal polynomials:
$$ \int_K  (y_n(z),y_n'(z),...,y_n^{(\rho)}(z)) M_0(z) 
\overline{
\left(
\begin{array}{cccc} y_m(z)\\
y_m'(z)\\
\vdots\\
y_m^{(\rho)}(z)
\end{array}
\right)
}
d\mu
= 
$$
\begin{equation}
\label{f2_10}
= A_n \delta_{n,m},\qquad               A_n>0,\quad           n,m\in\mathbb{Z}_+.
\end{equation}
Let us show that the system $\{ y_n(z) \}_{n=0}^\infty$
is related to some orthogonal systems in the direct sum of the scalar
$L^2_\mu$ spaces.

Suppose that the matrix function $M_0(z)$ admits the following factorization:
\begin{equation}
\label{f2_15}
M_0(z) = G(z) G^*(z),
\end{equation}
where $G(z) = (g_{l,k}(z))_{0\leq l\leq\rho,\ 0\leq k\leq \beta}$ is a measurable $\mathbb{C}_{\rho\times\beta}$-valued
function on $K$; $\beta\in\mathbb{N}$.
In the case $M_0 = M'_\tau$, one possible choice of $G(z)$ is given by the square root of $M'_\tau$, 
which is known to be measurable, see~\cite{cit_7000__Rosenberg_1964}.
However, as we shall see below it is better to choose $\beta$ as small as possible. 
For an arbitrary function $f(z)\in A^2(M)$, we denote
\begin{equation}
\label{f2_20}
g_{f;k}(z) := \sum_{l=0}^\rho g_{l,k}(z) f^{(l)}(z),\qquad k=0,1,...,\beta;\ z\in K; 
\end{equation}
\begin{equation}
\label{f2_25}
\vec g_{f}(z) := (g_{f;0}(z),g_{f;1}(z),...,g_{f;\beta}(z)),\qquad z\in K. 
\end{equation}
Set
\begin{equation}
\label{f2_27}
L^2_{\beta;\mu} = L^2_{\beta;\mu,K} := \bigoplus_{j=0}^\beta L^2_{\mu,K}.
\end{equation}
For arbitrary functions $f(z),u(z)\in A^2(M)$, we may write:
$$ (f,u)_{W^2(M)} = 
\int_K  (f(z),f'(z),...,f^{(\rho)}(z)) G(z) G^*(z) 
\overline{
\left(
\begin{array}{cccc} u(z)\\
u'(z)\\
\vdots\\
u^{(\rho)}(z)
\end{array}
\right)
}
d\mu
= 
$$
\begin{equation}
\label{f2_29}
=
\int_K \vec g_{f}(z) \left( \vec g_{u}(z) \right)^*  d\mu. 
\end{equation}
Taking $u(z) = f(z)$, we conclude that $g_{f;k}(z)\in L^2_\mu$, $k=0,1,...,\beta$.  
Then $\vec g_{f}(z)\in L^2_{\beta;\mu}$, and
\begin{equation}
\label{f2_33}
(f,u)_{W^2(M)} = 
(\vec g_{f}(z), \vec g_{u}(z))_{L^2_{\beta;\mu}},\qquad \forall f,u\in A^2(M). 
\end{equation}
In particular, for the Sobolev orthogonal polynomials $\{ y_n(z) \}_{n=0}^\infty$ we have the following property:
\begin{equation}
\label{f2_35}
A_n \delta_{n,m} = (y_n,y_m)_{W^2(M)} = 
(\vec g_n(z), \vec g_m(z))_{L^2_{\beta;\mu}},\qquad n,m\in \mathbb{Z}_+, 
\end{equation}
where  
\begin{equation}
\label{f2_38}
\vec g_j(z) := \vec g_{y_j}(z),\qquad j\in\mathbb{Z}_+. 
\end{equation}
Consequently, \textit{Sobolev orthogonal polynomials $\{ y_n(z) \}_{n=0}^\infty$ lead to an orthogonal system of
functions $\{ \vec g_n(z) \}_{n=0}^\infty$ in the space $L^2_{\beta;\mu}$}.
This property can be used to construct Sobolev orthogonal polynomials using known orthogonal systems of functions in
$L^2_{\beta;\mu}$. We come to the following question.

\noindent
\textbf{Question 1.} \textit{ Let a space $L^2(M)$ on $K\in\mathfrak{B}(\mathbb{C})$ be given. Suppose that the measurable factorization~(\ref{f2_15})
with some $\beta\in\mathbb{N}$ holds.
Fix some orthogonal system of
functions $\{ \vec h_n(z) \}_{n=0}^\infty$, $\vec h_n(z) = (h_{n;0}(z),...,h_{n;\beta}(z))$, in the space $L^2_{\beta;\mu}$.
Do there exist complex polynomials $p_n(z)$, $\deg p_n = n$, $n=0,1,2,...$, which are solutions of the following system of
differential equations: 
\begin{equation}
\label{f2_45}
h_{n;k}(z) = \sum_{l=0}^\rho g_{l,k}(z) p_n^{(l)}(z),\qquad k=0,1,...,\beta;\ z\in K. 
\end{equation}
}
If an answer on Question~1 is affirmative, then $\{ p_n(z) \}_{n=0}^\infty$ belong to $A^2(M)$.
In fact, this follows from the existence of integrals on the right of~(\ref{f2_29}) for $\vec g_n(z) = \vec h_n(z)$.
One can repeat the constructions after~(\ref{f2_15}), taking into account~(\ref{f2_45}), to conclude that 
$\{ p_n(z) \}_{n=0}^\infty$ are Sobolev orthogonal polynomials.

Differential equations~(\ref{f2_45}) can be used in the search for an explicit representation of $p_n$.
In the case when $\beta = 0$, $g_{l,k}(z), h_{n;k}(z) \in\mathbb{P}$, the following proposition is useful.

\begin{proposition}
\label{p2_1}
Let $D$ be a linear differential operator of order $r\in\mathbb{N}$, with complex polynomial coefficients:
\begin{equation}
\label{f2_47}
D = \sum_{k=0}^r d_k(z) \frac{ d^k }{ dz^k},\quad d_k(z)\in\mathbb{P}.
\end{equation}
Let $\{ u_n(z) \}_{n=0}^\infty$, $\deg u_n = n$, be an arbitrary set of complex polynomials.
The following statements are equivalent:

\begin{itemize}
\item[(A)] The following equation:
\begin{equation}
\label{f2_50}
D y(z) = u_n(z),
\end{equation}
for each $n\in\mathbb{Z}_+$, has a complex polynomial solution $y(z)=y_n(z)$ of degree $n$;

\item[(B)] $D z^n$ is a complex polynomial of degree $n$, $\forall n\in\mathbb{Z}_+$; 

\item[(C)] 
The following conditions hold:
\begin{equation}
\label{f2_53}
\deg d_k \leq k,\qquad 0\leq k\leq r;
\end{equation}
\begin{equation}
\label{f2_55}
\sum_{j=0}^r [n]_j d_{j,j}\not= 0,\qquad n\in\mathbb{Z}_+,
\end{equation}
where $d_{j,l}$ means the coefficient by $z^l$ of the polynomial $d_j$.
\end{itemize}

If one of the statements~$(A),(B),(C)$ holds true, then for each $n\in\mathbb{Z}_+$, the solution of~(\ref{f2_50}) is unique.

\end{proposition}
\textbf{Proof.} 
$(A)\Rightarrow(B)$.
For each $n\in\mathbb{Z}_+$, the polynomial $z^n$ can be expanded as a linear combination of $y_j(z)$, $0\leq j\leq n$:
\begin{equation}
\label{f2_56}
z^n = \sum_{j=0}^n \varphi_{n,j} y_j(z),\qquad  \varphi_{n,j}\in\mathbb{C},\ \varphi_{n,n}\not=0.
\end{equation}
Then
$$ D z^n = \sum_{j=0}^n \varphi_{n,j} D y_j(z) = \sum_{j=0}^n \varphi_{n,j} u_j(z), $$
is a polynomial of degree $n$.

\noindent
$(B)\Rightarrow(A)$.
Set
$$ t_n(z) := D z^n,\qquad n\in\mathbb{Z}_+. $$
By condition~$(B)$ we see that $\deg t_n = n$, $n\in\mathbb{Z}_+$.
We can expand the polynomial $u_n(z)$ as a linear combination of $t_j(z)$, $0\leq j\leq n$:
\begin{equation}
\label{f2_57}
u_n(z) = \sum_{j=0}^n \xi_{n,j} t_j(z),\qquad  \xi_{n,j}\in\mathbb{C},\ \xi_{n,n}\not=0;\quad n\in\mathbb{Z}_+.
\end{equation}
Set
\begin{equation}
\label{f2_59}
y_n(z) = \sum_{j=0}^n \xi_{n,j} z^j,\qquad  n\in\mathbb{Z}_+.
\end{equation}
Then $\deg y_n = n$, and
$$ D y_n(z) = \sum_{j=0}^n \xi_{n,j} D z^j = \sum_{j=0}^n \xi_{n,j} t_j(z) = u_n(z), $$
for all $n\in\mathbb{Z}_+$.

\noindent
$(B)\Rightarrow(C)$.
We shall check condition~(\ref{f2_53}) by the induction argument.
For $k=0$, we have: 
$$ D 1 = d_0(z), $$
and $D 1$ has degree $0$ by condition $(B)$. Thus, $\deg d_0 = 0$.
Suppose that condition~(\ref{f2_53}) holds for $k\in\mathbb{Z}_{0,l}$, with some $l\in\mathbb{Z}_{0,r-1}$.
Let us verify condition~(\ref{f2_53}) for $k=l+1$. We may write:
\begin{equation}
\label{f2_62}
D z^{l+1} = \sum_{k=0}^{l} d_k(z) \frac{ d^k }{ dz^k} z^{l+1} + d_{l+1}(z) (l+1)!. 
\end{equation}
By the induction assumptions the sum $\sum_{k=0}^{l} ...$, on the right has degree $\leq l+1$.
The degree of the left-hand side is equal to $l+1$ by condition~$(B)$.
Therefore $\deg d_{l+1} \leq l+1$. Consequently, condition~(\ref{f2_53}) holds.

For each $n\in\mathbb{Z}_+$, we may write:
\begin{equation}
\label{f2_64}
D z^n = \sum_{k=0}^r d_k(z) \frac{ d^k }{ dz^k } z^n = \sum_{k=0}^r d_k(z) [n]_k z^{n-k}. 
\end{equation}
By condition~(\ref{f2_53}) we see that the polynomial on the right is of degree $\leq n$.
The coefficient of $z^n$ is equal to $\sum_{k=0}^r d_{k,k} [n]_k$.
By condition~$(B)$ the left-hand side of~(\ref{f2_64}) has degree $n$. Therefore
relation~(\ref{f2_55}) holds.

\noindent
$(C)\Rightarrow(B)$. By~(\ref{f2_64}) and conditions~(\ref{f2_53}),(\ref{f2_55}) we see that $Dz^n$ has degree $n$.

Let us check the last statement of the theorem.
Suppose to the contrary that for some $n_0\in\mathbb{Z}_+$, there exists another polynomial solution $v_{n_0}(z)$ ($\deg v_{n_0} = n_0$) of~(\ref{f2_50}),
which is different from $y_{n_0}(z)$.
Then
\begin{equation}
\label{f2_65}
D( y_{n_0}(z) - v_{n_0}(z) ) = D y_{n_0}(z) - D v_{n_0}(z) = u_{n_0}(z) - u_{n_0}(z) = 0. 
\end{equation}
Observe that $y_{n_0}(z) - v_{n_0}(z)$, is a non-zero polynomial of degree $k_0$, $0\leq k_0\leq n_0$. 
By condition~$(B)$ we conclude that $D(y_{n_0}(z) - v_{n_0}(z))$ should have degree $k_0$ as well.
This contradicts to relation~(\ref{f2_65}).
Thus, $y_n(z)$ is a unique solution of~(\ref{f2_50}).
$\Box$

Observe that condition~(\ref{f2_55}) holds true, if the following simple condition holds:
\begin{equation}
\label{f2_67}
d_{0,0}>0,\quad d_{j,j}\geq 0,\qquad j\in\mathbb{Z}_{1,r}.
\end{equation}
Thus, there exists a big variety of linear differential operators with polynomial coefficients which
have the property~$(A)$.

Let us show how Proposition~\ref{p2_1} can be applied to a question of the density of polynomials in $W^2(M)$.
Let $r\in\mathbb{N}$, be an arbitrary number, and
\begin{equation}
\label{f2_70}
G(z) = 
\left(
\begin{array}{cccc} d_0(z)\\
d_1(z)\\
\vdots\\
d_r(z) \end{array}
\right),
\end{equation}
where $d_j(z)\in\mathbb{P}$, $j\in\mathbb{Z}_{0,r}$, satisfy condition~$(C)$ of Proposition~\ref{p2_1}.
Set
\begin{equation}
\label{f2_72}
M_0(z) = G(z) G^*(z).
\end{equation}
Let $K$ be an arbitrary Borel subset of the complex plane. Let $\mu$ be an arbitrary (non-negative) finite measure on $\mathfrak{B}(K)$,
having all finite power moments on $K$: 
\begin{equation}
\label{f2_73}
\int_K z^j d\mu < \infty,\qquad j\in\mathbb{Z}_+. 
\end{equation}
Define the following matrix measure:
\begin{equation}
\label{f2_75}
M(\delta) = \int_\delta M_0(z) d\mu,\qquad \delta\in\mathfrak{B}(K).
\end{equation}
Consider the corresponding spaces $L^2(M)$ and $W^2(M)$. 
Observe that \textit{we do not assume the validity of Condition~(\ref{f2_7})}.
We may repeat our constructions after~(\ref{f2_10}) up to~(\ref{f2_33}), except of those concerning Sobolev orthogonal polynomials.

Suppose that complex polynomials are dense in $L^2_{\mu}$. 
Choose an arbitrary $f\in A^2(M)$.
For an arbitrary $\varepsilon > 0$, there exists
$p_\varepsilon\in\mathbb{P}$, such that 
\begin{equation}
\label{f2_79}
\| \vec g_{f}(z) - p_\varepsilon(z) \|_{L^2_{\mu}} < \varepsilon. 
\end{equation}
By~condition~$(A)$ of Proposition~\ref{p2_1} there exists $w\in\mathbb{P}$, such that
$$ \vec g_{w}(z) = p_\varepsilon(z). $$
By~(\ref{f2_33}) we may write
$$ \| f-w \|_{W^2(M)}^2 = 
(\vec g_{f-w}(z), \vec g_{f-w}(z))_{L^2_{\mu}} = $$ 
\begin{equation}
\label{f2_82}
= (\vec g_{f}(z) - \vec g_{w}(z), \vec g_{f}(z) - \vec g_{w}(z))_{L^2_{\mu}} =
\| \vec g_{f}(z) - p_\varepsilon(z) \|_{L^2_{\mu}}^2 < \varepsilon^2.
\end{equation}
Thus, polynomials are dense in $W^2(M)$.

\begin{corollary}
\label{c2_1}
Let $K$ be an arbitrary Borel subset of the complex plane. Let $\mu$ be an arbitrary (non-negative) finite measure on $\mathfrak{B}(K)$,
having all finite power moments~(\ref{f2_73}).
Let $r\in\mathbb{N}$, and $G(z)$ be defined by~(\ref{f2_70}),
where $d_j(z)\in\mathbb{P}$, $j\in\mathbb{Z}_{0,r}$, satisfy condition~$(C)$ of Proposition~\ref{p2_1}.
Define a matrix measure $M$ by relations~(\ref{f2_72}),(\ref{f2_75}).
If complex polynomials are dense in $L^2_{\mu}$, then complex polynomials are dense in $W^2(M)$.
In particular, if $K = \mathbb{T}$, and $\mu$ is nontrivial, then the Szeg\"o condition:
\begin{equation}
\label{f2_85}
\int_0^{2\pi} \log w(\theta) \frac{d\theta}{2\pi} = -\infty,
\end{equation}
is sufficient for the density of polynomials in the corresponding space $W^2(M)$.
Here $w = 2\pi d\mu_{ac}/d\theta$.

\end{corollary}
\textbf{Proof.} The proof follows from the preceding considerations.
$\Box$

In 1992 Klotz introduced $L^p(M)$ spaces, for $0<p<\infty$, having square $(\rho+1)\times(\rho+1)$ 
matrix functions as representatives of elements ($\rho\in\mathbb{Z}_+$), 
see~\cite{cit_1700__Klotz_1992}.
Denote by $A^p(M)$ a set in $L^2(M)$ including those classes of the equivalence $[\cdot]$ which possess 
a representative of the following form: 
\begin{equation}
\label{f2_90}
F(z) = \left(
\begin{array}{cccc}
f(z) & f'(z) & \cdots & f^{(\rho)}(z) \\
0 & 0 & \cdots & 0 \\
\vdots & \vdots & \ddots & \vdots \\
0 & 0 & \cdots & 0 \end{array}
\right).
\end{equation}
By $W^p(M)$ we denote the closure of $A^p(M)$ in the norm of $L^p(M)$. 
The subspace $W^p(M)$ is said to be \textbf{the Sobolev space of index $p$ with the matrix measure $M$}.
It is of interest to study the approximation by polynomials in spaces $W^p(M)$.
It is possible that ideas of Xu from~\cite{cit_9200__Xu_2018} can be applied in this case.


Let us now show how to use measurable factorizations of type~(\ref{f2_15}) in order to obtain new measures of orthogonality
for Sobolev orthogonal polynomials.
Let $M(\delta)$ be a non-negative Hermitian-valued measure on $(K,\mathfrak{B}(K))$, 
with an arbitrary Borel subset $K$ of the complex plane and an arbitrary $\rho\in\mathbb{N}$. 
Let $\mu$ be an admissible measure, $M_0(z)$ be the corresponding matrix weight, and $1,z,z^2,...$, all belong to $W^2(M)$. 
As before, we suppose that condition~(\ref{f2_7}) holds. 
Let
$\{ y_n(z) \}_{n=0}^\infty$, $\deg y_n=n$, be the corresponding sequence of Sobolev orthogonal polynomials~(see relation~(\ref{f2_10})).
Suppose that some factorization~(\ref{f2_15}) holds, and, in notations~(\ref{f2_20}),(\ref{f2_25}), we have relation~(\ref{f2_33}).
We assume that the following condition is valid:

\noindent
\textbf{Condition 1.} \textit{For some $k_0\in\mathbb{Z}_{0,\beta}$, the functions $g_{l,k_0}(z)$, $0\leq l\leq \rho$, all have $\mu$-measurable
derivatives. Functions $\{ g_{y_n;k_0}'(z) \}_{n=0}^\infty$ are pairwise orthogonal:
\begin{equation}
\label{f2_92}
\left(
g_{y_n;k_0}'(z), g_{y_m;k_0}'(z)
\right)_{L^2_{\mu,K}} = B_n\delta_{n,m},\qquad B_n\geq 0,\quad n,m\in\mathbb{Z}_+.
\end{equation}
}

Notice that $B_n$ in~(\ref{f2_92}) can be zero.
Differentiating relation~(\ref{f2_20}) (with $k=k_0$, $f=y_n$, $n\in\mathbb{Z}_+$) we obtain that
\begin{equation}
\label{f2_94}
g_{y_n;k_0}'(z) = \sum_{l=0}^\rho d_l(z) y_n^{(l)}(z) + g_{\rho,k_0}(z) y_n^{(\rho+1)}(z),
\end{equation}
where
\begin{equation}
\label{f2_96}
d_l(z) := g_{l,k_0}'(z) + g_{l-1,k_0}(z),\quad l\in\mathbb{Z}_{0,\rho},\quad g_{-1,k_0}(z) := 0.
\end{equation}
Set
\begin{equation}
\label{f2_97}
\widetilde G(z) = (\widetilde g_{l,k}(z))_{0\leq l\leq\rho+1,\ 0\leq k\leq \beta+1},
\end{equation}
where
\begin{equation}
\label{f2_98}
\widetilde g_{l,k}(z) =
\left\{ 
\begin{array}{cccc}
g_{l,k}(z), & 0\leq l\leq\rho,\ 0\leq k\leq \beta\\
0, & l=\rho+1,\ 0\leq k\leq \beta\\
d_l(z), & 0\leq l\leq\rho,\ k = \beta+1\\
g_{\rho,k_0}(z), & l=\rho+1,\ k=\beta+1
\end{array}
\right..
\end{equation}
Then
$$ \int_K  (y_n(z),y_n'(z),...,y_n^{(\rho)}(z),y_n^{(\rho+1)}(z)) \widetilde G(z) \widetilde G^*(z) 
\overline{
\left(
\begin{array}{ccccc} y_m(z)\\
y_m'(z)\\
\vdots\\
y_m^{(\rho)}(z)\\
y_m^{(\rho+1)}(z)
\end{array}
\right)
}
d\mu
= $$
$$ = \int_K  \left( \vec g_{y_n}(z), g_{y_n;k_0}'(z) \right)
\left( \vec g_{y_m}(z), g_{y_m;k_0}'(z) \right)^*
d\mu
= $$
$$ = \int_K  \vec g_{y_n}(z) (\vec g_{y_m}(z))^*
d\mu
+
\int_K  g_{y_n;k_0}'(z)
\overline{ g_{y_m;k_0}'(z) }
d\mu
= 
$$
\begin{equation}
\label{f2_100}
= (A_n+B_n)\delta_{n,m},\quad n,m\in\mathbb{Z}_+.
\end{equation}
Thus, we have obtained new orthogonality relations for $y_n$. The following example illustrates this construction.

\begin{example}
\label{e2_1}
Consider polynomials from Theorem~2.1 in~\cite{cit_9000__Zagorodnyuk_CMA_2020}:
$$ y_n(z) = y_n(1;z) = -\frac{1}{n!} z^n {}_2 F_0\left(
-n,1;-;-\frac{1}{z}
\right),\qquad z\in\mathbb{C}\backslash\{ 0 \},\ n\in\mathbb{Z}_+.
$$
They satisfy the following orthogonality relations:
$$ \int_\mathbb{T}  (y_n(z),y_n'(z)) \mathcal M
\overline{
\left(
\begin{array}{cc} y_m(z)\\
y_m'(z)\end{array}
\right)
}
d\mu_0
= \delta_{n,m},\qquad n,m\in\mathbb{Z}_+, $$
where
$\mathcal M = \left(
\begin{array}{cc}
1 & -1\\
-1 & 1\end{array}
\right)$, and $\mu_0$ is the Lebesgue measure on $[0,2\pi)$ (with the identification $z=e^{i\theta}$).
Observe that
$$ \mathcal M = 
\left(
\begin{array}{cc}
-1\\
1\end{array}
\right) (-1, 1). $$
Set $\rho=1$,
$G(z) = (-1,1)^T$, and repeat constructions~(\ref{f2_70})--(\ref{f2_75}) to define a matrix measure $M$ on $\mathfrak{B}(\mathbb{T})$.
All necessary assumptions for the above construction of a new measure, including Condition~1, are satisfied.
In fact, we have $K=\mathbb{T}$, $M_0(z)=\mathcal M$, $\mu$ is $\mu_0$, restricted to $\mathfrak{B}(\mathbb{T})$.
Relation~(\ref{f2_20}) takes the following form:
$$ g_{y_n;0}(z) = -y_n(z) + y_n'(z) = z^n,\qquad n\in\mathbb{Z}_+, $$
where the last equality is from~\cite[formula (2.13)]{cit_9000__Zagorodnyuk_CMA_2020}.
Then
$$ \widetilde G(z) =
\left(
\begin{array}{ccc}
-1 & 0\\
1 & -1\\
0 & 1\end{array}
\right),\quad 
\widetilde G(z) \widetilde G(z)^* =
\left(
\begin{array}{ccc}
1 & -1 & 0\\
-1 & 2 & -1\\
0 & -1 & 1\end{array}
\right). $$
Finally, we have
$B_n = 
\left\{
\begin{array}{cc}
n^2, & if\ n\geq 1\\
0, & if\ n=0\end{array}
\right.$; $A_n=1$, $n\in\mathbb{Z}_+$.

\end{example}

\section{Collections of Sobolev orthogonal polynomials with explicit integral representations.}


Suppose that conditions $(A),(B),(C)$ of Proposition~\ref{p2_1} are satisfied for a differential operator $D$,
and a set of polynomials $\{ u_n(z) \}_{n=0}^\infty$. Let us study equation~(\ref{f2_50}).

In order to obtain an explicit representation for $y_n$, it is convenient to equate the corresponding powers in~(\ref{f2_50})
and solve the corresponding linear system of equations for the unknown coefficients of $y_n$.
Such an idea was used in~\cite{cit_5_Azad}. However, as we shall see, this way leads to computational difficulties.
Fix an arbitrary $n\in\mathbb{Z}_+$. Let
\begin{equation}
\label{f3_10}
y(z) = \sum_{j=0}^n \mu_{n,j} z^j,\qquad \mu_{n,j}\in\mathbb{C}.
\end{equation}
Then
$$ D y(z) = \sum_{j=0}^n \mu_{n,j} D z^j = \sum_{j=0}^n \sum_{l=0}^r \mu_{n,j} d_l(z) [j]_l z^{j-l} = $$
$$ = \sum_{j=0}^n \sum_{l=0}^r \sum_{k=0}^l \mu_{n,j} [j]_l d_{l,k} z^{j-l+k}. $$
Using the change of indices in the last sum: $t= j-l+k$, $k=l+t-j$, we get
\begin{equation}
\label{f3_15}
D y(z) = \sum_{j=0}^n \sum_{l=0}^r \sum_{t=j-l}^j \mu_{n,j} [j]_l d_{l,l+t-j} z^{t} = 
\sum_{j=0}^n \sum_{l=0}^r \sum_{t=0}^j \mu_{n,j} [j]_l d_{l,l+t-j} z^{t}.
\end{equation}
In the last sum we replaced the lower bound $t=j-l$, by $t=0$. In the case $j\leq l$, this is correct, since $Dy$ is a polynomial.
If $j>l$, then there appear new terms with $0\leq t\leq j-l-1$.
For these terms we have $l+t-j\leq -1 < 0$. Thus, we define $d_{l,m}:=0$, for $l\in\mathbb{Z}_{0,r}$, $m\in\mathbb{Z}:\ m<0$. 

Changing the order of summation for $j$ and $t$ we obtain that
\begin{equation}
\label{f3_17}
D y(z) = 
\sum_{t=0}^n z^t \left(
\sum_{l=0}^r \sum_{j=t}^n \mu_{n,j} [j]_l d_{l,l+t-j}
\right).
\end{equation}
Let
\begin{equation}
\label{f3_19}
u_n(z) = \sum_{j=0}^n a_{n,j} z^j,\qquad a_{n,j}\in\mathbb{C},\ a_{n,n}\not=0.
\end{equation}
Equation~(\ref{f2_50}) is now equivalent to the following system of linear algebraic equations:
\begin{equation}
\label{f3_22}
\mu_{n,n} \sum_{l=0}^r [n]_l d_{l,l} = a_{n,n},
\end{equation}
$$ \mu_{n,n-k} \sum_{l=0}^r [n-k]_l d_{l,l} + \sum_{j=n-k+1}^n \mu_{n,j} \sum_{l=0}^r [j]_l d_{l,l+n-k-j} = a_{n,n-k}, $$
\begin{equation}
\label{f3_25}
k=1,2,...,n.
\end{equation}
It is clear that one can find $\mu_{n,n},\mu_{n,n-1},...,\mu_{n,0}$, step by step. However, for large $n$ this leads to 
huge expressions. The use of Cramer's rule leads to large determinants.
Thus, it is important to point out such cases for $D$, which admit convenient explicit representations for $y_n$.
As we shall show, this can be successfully done, if the operator $D$ has constant coefficients.

Fix arbitrary $r\in\mathbb{N}$, and $\alpha\in\mathbb{R}$.
We shall begin with the following simple differential operator:
\begin{equation}
\label{f3_30}
D = \alpha \frac{ d^r }{ dx^r } + 1.
\end{equation}
In this case the system~(\ref{f3_22}),(\ref{f3_25}) can be solved directly (as a two-term recurrence relation).
The properties of the corresponding polynomial solutions will give a key to the announced general case of constant coefficients.

Observe that $D$ satisfies condition~(\ref{f2_67}). Therefore it satisfies conditions $(A),(B),(C)$ of Proposition~\ref{p2_1}.
At first, we shall consider the case: $u_n(z)=z^n$, $n\in\mathbb{Z}_+$.
Equation~(\ref{f3_22}) implies that $\mu_{n,n}=1$, while equations~(\ref{f3_25}) take the following form:
\begin{equation}
\label{f3_35}
\mu_{n,n-k} + \sum_{j=n-k+1}^n \mu_{n,j} \sum_{l=0}^r [j]_l d_{l,l+n-k-j} = 0,\qquad k=1,2,...,n.
\end{equation}
Using the change of index: $s=n-k$, $k=n-s$, we may rewrite~(\ref{f3_35}) in the following form:
\begin{equation}
\label{f3_37}
\mu_{n,s} + \sum_{j=s+1}^n \mu_{n,j} \sum_{l=0}^r [j]_l d_{l,l-j+s} = 0,\qquad s=0,1,...,n-1.
\end{equation}
If $l\notin\{ 0,r \}$, then $d_{l,l-j+s}=0$. If $l=0$, then $d_{l,l-j+s}=0$, since $s-j < 0$.
Finally, if $l=r$, then
\begin{equation}
\label{f3_45}
d_{l,l-j+s}= d_{r,r-j+s}= \left\{
\begin{array}{cc}
\alpha, & \mbox{if $j=s+r$}\\
0, & \mbox{if $j\not=s+r$}\end{array}
\right..
\end{equation}
Therefore
\begin{equation}
\label{f3_48}
\left\{
\begin{array}{cc}
\mu_{n,s} + \mu_{n,s+r} \alpha [s+r]_r = 0, & \mbox{if $s\leq n-r$}\\
\mu_{n,s} = 0, & \mbox{if $s>n-r$}\end{array}
\right.;\qquad s=0,1,...,n-1.
\end{equation}
By the induction argument one can check that
\begin{equation}
\label{f3_50}
\mu_{n,n-kr} = (-\alpha)^k [n]_r [n-r]_r ... [n-(k-1)r]_r,\quad k\in\mathbb{N}:\ n-kr\geq 0. 
\end{equation}
If $r=1$, we have calculated all the coefficients of $y(z)$.
If $r>1$, then by~(\ref{f3_48}) we see that 
$$ \mu_{n,n-r+1} = ... = \mu_{n,n-1} = 0. $$
Using the first equation in~(\ref{f3_48}) and the induction we conclude that the rest of coefficients of $y$ are zeros.

\begin{theorem}
\label{t3_1}
Let $r\in\mathbb{N}$, $\alpha\in\mathbb{R}$, be arbitrary numbers.
Polynomials
\begin{equation}
\label{f3_54}
y_n(z) = y_n(r,\alpha;z) := z^n + n! \sum_{k=1}^{\left[ \frac{n}{r} \right]} (-\alpha)^k \frac{ z^{n-kr} }{ (n-kr)! },\quad n\in\mathbb{Z}_+,
\end{equation}
have the following properties:
\begin{itemize}
\item[(i)] They satisfy the following differential equation:
\begin{equation}
\label{f3_57}
\alpha z y_n^{(r+1)}(z) + z y_n'(z) = n \left(
\alpha y_n^{(r)}(z) + y_n(z)
\right),\qquad n\in\mathbb{Z}_+;
\end{equation}

\item[(ii)] Polynomials $y_n$ obey the following mixed relation:
\begin{equation}
\label{f3_59}
\alpha y_{n+1}^{(r)}(z) +  y_{n+1}(z) = z \left(
\alpha y_n^{(r)}(z) + y_n(z)
\right),\qquad n\in\mathbb{Z}_+;
\end{equation}

\item[(iii)] Polynomials $y_n$ are Sobolev orthogonal polynomials on the unit circle:
\begin{equation}
\label{f3_61}
\int_{\mathbb{T}} \left( y_n(z), y_n'(z),..., y_n^{(r)}(z) \right) M \overline{
\left( \begin{array}{cccc} y_m(z) \\
y_m'(z) \\
\vdots \\
y_m^{(r)}(z) \end{array} \right)
} 
\frac{1}{2\pi} d\theta = \delta_{n,m},\qquad n,m\in\mathbb{Z}_+,
\end{equation}
where
\begin{equation}
\label{f3_62}
M = (1,0,...,0,\alpha)^T (1,0,...,0,\alpha).
\end{equation}
Here $d\theta$ means the Lebesgue measure on $[0,2\pi)$ (which may be identified with $\mathbb{T}$, by $z=e^{i\theta}$).

\item[(iv)] If $\alpha\not=0$; $l\in\mathbb{Z}_{0,r-1}$, then
\begin{equation}
\label{f3_63}
\frac{ \alpha_r^{rm+l} }{ (rm+l)! } y_{rm+l}(r,\alpha;z) \rightarrow_{m\rightarrow+\infty}
\frac{1}{r} \sum_{k=0}^{r-1} \varepsilon^{-lk} e^{ \alpha_r \varepsilon^k z},\qquad \forall z\in\mathbb{C},
\end{equation}
where $\varepsilon$ is a primitive $r$-th root of unity; $\alpha_r =
\left(
-\frac{1}{\alpha}
\right)^{ \frac{1}{r} }$, with an arbitrary fixed value of the $r$-th root.
The convergence in~(\ref{f3_63}) is uniform on compact subsets of the complex plane.

\item[(v)] If $\alpha\leq -1$; $m\in\mathbb{Z}_+$, then 
$y_{rm}(r,\alpha;z)$ has all its roots in $\{ z\in\mathbb{C}:\ |z|\geq 1 \}$.
If $r\geq 2$, $\alpha\leq -1$; $m\in\mathbb{Z}_+$, $l\in\mathbb{Z}_{1,r-1}$,
then 
$y_{rm}(r,\alpha;z)$ has a root $z=0$ of multiplicity $l$, and all other roots in $\{ z\in\mathbb{C}:\ |z|\geq 1 \}$.
\end{itemize}

\end{theorem}
\textbf{Proof.} 
Observe that polynomials $y_n(z)$ from~(\ref{f3_54}) are those polynomials which we have constructed before
the statement of the theorem (after some simplifications for $\mu_{n,n-kr}$).
Using relation~(\ref{f2_50}) (with $y=y_n$, $u_n(z)=z^n$) and relation~(\ref{f1_50}) we obtain
property~$(i)$ of the theorem.
By~(\ref{f2_50}) and 
$$ z^{n+1} = z z^n,\qquad n\in\mathbb{Z}_+, $$
we conclude that property~$(ii)$ holds.
By~(\ref{f2_50}) and the orthogonality relations for $z^n$ we obtain property~$(iii)$.

Let us check relation~(\ref{f3_63}). By~(\ref{f3_54}) we may write
\begin{equation}
\label{f3_63_2}
\frac{1}{(rm+l)!} y_{rm+l}(r,\alpha;z) = \sum_{j=0}^m (-\alpha)^{m-j} \frac{ z^{rj+l} }{ (rj+l)! },\quad
m\in\mathbb{Z}_+,\ l\in\mathbb{Z}_{0,r-1}. 
\end{equation}
If $\alpha\not=0$, then
$$ \alpha_r^l \frac{1}{ (-\alpha)^m (rm+l)! }
y_{rm+l}(r,\alpha;z) = \sum_{j=0}^m \frac{ (\alpha_r z)^{rj+l} }{ (rj+l)! }
\rightarrow_{m\rightarrow+\infty} $$
\begin{equation}
\label{f3_64}
\rightarrow_{m\rightarrow+\infty} \sum_{j=0}^\infty \frac{ (\alpha_r z)^{rj+l} }{ (rj+l)! },\quad \forall z\in\mathbb{C}.
\end{equation}
Notice that
\begin{equation}
\label{f3_65}
\sum_{j=0}^\infty \frac{ (\alpha_r z)^{rj+l} }{ (rj+l)! } =
\sum_{n=0}^\infty a_n \alpha_r^n \frac{ z^n }{ n! }, 
\end{equation}
where $a_n=1$, if $n=l+tr$, $t\in\mathbb{Z}_+$, and $a_n=0$, otherwise.
By Abel's theorem for power series it follows that the convergence is uniform on compact subsets of $\mathbb{C}$.

\noindent
We shall use the idea of operators $T_{m,N}$ from~\cite[p. 88]{cit_100_Duran}. Observe that
\begin{equation}
\label{f3_65_1}
\frac{1}{r} \sum_{k=0}^{r-1} \varepsilon^{-lk} e^{ \alpha_r z \varepsilon^k} =
\sum_{n=0}^\infty a_n \frac{(\alpha_r z)^n}{n!}, 
\end{equation}
where we applied the following relation:
$$ \sum_{k=0}^{r-1} \varepsilon^{k(n-l)} = \left\{
\begin{array}{cc}
r, & if\ N=rj+l\\
0, & otherwise
\end{array}
\right.. $$
Comparing~(\ref{f3_65}) and (\ref{f3_65_1}) we finish the proof of~$(iv)$.

Let us check property~$(v)$. By~(\ref{f3_63_2}) we may write
$$ \frac{1}{(rm+l)!} y_{rm+l}(r,\alpha;z) = z^l f(z^{r}),\quad $$
$$ f(w):= \sum_{j=0}^m (-\alpha)^{m-j} \frac{ w^j }{ (rj+l)! };\ m\in\mathbb{Z}_+,\ l\in\mathbb{Z}_{0,r-1}. $$
We can apply the Enestr\"om--Kakeya Theorem~(\cite[p. 136]{cit_8100_Marden}) for the polynomial $f(w)$
to obtain the required assertions.
$\Box$

Polynomials $y_n(z)=y_n(r,\alpha;z)$ from~(\ref{f3_54}) allow to construct Sobolev orthogonal polynomials from any sequence of
orthogonal polynomials on the unit circle or on the real line.

\begin{corollary}
\label{c3_1}
Let $\mathcal{K}$ denote the real line or the unit circle, and $r\in\mathbb{N}$, $\alpha\in\mathbb{R}$.
Let $\mu$ be a (non-negative) measure on $\mathcal{K}$ (we
assume that it is defined at least on $\mathfrak{B}(\mathcal{K})$). 
Denote by
$p_n$ orthogonal polynomials on $\mathcal{K}$ with respect to $\mu$ (the positivity of leading coefficients
is not assumed):
\begin{equation}
\label{f3_67}
\int_{\mathcal{K}} p_n(z) \overline{p_m(z)} d\mu = A_n \delta_{n,m},\qquad A_n>0,\ n,m\in\mathbb{Z}_+.
\end{equation}
Let
\begin{equation}
\label{f3_73_1}
p_n(z) = \sum_{j=0}^n \xi_{n,j} z^j,\qquad \xi_{n,j}\in\mathbb{C},\ \xi_{n,n}\not=0;\ n\in\mathbb{Z}_+.
\end{equation}
Polynomials
\begin{equation}
\label{f3_73_2}
\widehat y_n(z) = \widehat y_n(r,\alpha;z) = \sum_{j=0}^n \xi_{n,j} 
y_j(r,\alpha;z),\qquad n\in\mathbb{Z}_+,
\end{equation}
are Sobolev orthogonal polynomials on $\mathcal{K}$:
\begin{equation}
\label{f3_75}
\int_{\mathcal{K}} \left( \widehat y_n(z), \widehat y_n'(z),..., \widehat y_n^{(r)}(z) \right) M \overline{
\left( \begin{array}{cccc} \widehat y_m(z) \\
\widehat y_m'(z) \\
\vdots \\
\widehat y_m^{(r)}(z) \end{array} \right)
} 
d\mu = A_n \delta_{n,m},\qquad n,m\in\mathbb{Z}_+,
\end{equation}
where $M$ is from~(\ref{f3_62}).

\end{corollary}
\textbf{Proof.}
Consider the operator $D$ from~(\ref{f3_30}).
Since $D$ is a linear operator on polynomials, then we may write:
$$ D \sum_{j=0}^n \xi_{n,j} y_j(r,\alpha;z) =
\sum_{j=0}^n \xi_{n,j} D y_j(r,\alpha;z) = \sum_{j=0}^n \xi_{n,j} z^j = p_n(z). $$
Thus, $\widehat y_n$ are solutions to equation~(\ref{f2_50}) with $u_n=p_n$, and $D$ from~(\ref{f3_30}).
Substitute for $p_n$ into the orthogonality relations~(\ref{f3_67}) to obtain relation~(\ref{f3_75}).
$\Box$

Notice that polynomials $y_n(1,-1;z)$ are close to a particular case of polynomials which appeared
in~\cite{cit_9000__Zagorodnyuk_CMA_2020} (where recurrence relations were established).
We shall now consider the special case $r=2$, $\alpha<0$. Our aim is to obtain some recurrence relations for $y_n(2,\alpha;z)$ in
this case. For convenience we denote
\begin{equation}
\label{f3_77}
w_n(z) = w_n(\alpha;z) := y_n(2,\alpha;z),\qquad n\in\mathbb{Z}_+;\ \alpha<0.
\end{equation}

\begin{theorem}
\label{t3_2}
Let $\alpha<0$, be an arbitrary number, and polynomials $w_n(z) = w_n(\alpha;z)$, $n\in\mathbb{Z}_+$, be defined by~(\ref{f3_77}).
The following statements hold:
\begin{itemize}
\item[(a)]  Polynomials $w_n$ have the following integral representation:
$$ w_n(t) = 
\frac{\beta}{2} e^{\beta t} 
\int_t^{+\infty} x^n e^{-\beta x} dx  
+\frac{\beta}{2} e^{-\beta t}
\int_{-\infty}^{t} x^n e^{\beta x} dx, $$ 
\begin{equation}
\label{f3_78}
t\in\mathbb{R},\ n\in\mathbb{Z}_+;\quad \beta:= \sqrt{ -\frac{1}{\alpha} };  
\end{equation}

\item[(b)] Polynomials $w_n$ satisfy the following recurrence relation:
$$ w_{n+1}(z) + \alpha n(n+1) w_{n-1}(z) = 
z
\left(
w_{n}(z) + \alpha (n-1)n w_{n-2}(z)
\right),\qquad $$
\begin{equation}
\label{f3_79}
n\in\mathbb{Z}_+,
\end{equation}
where $w_{-1}(z) = w_{-2}(z) = 0$.

\item[(c)] Polynomials $w_n$ have the following generating function:
\begin{equation}
\label{f3_80}
\sum_{n=0}^\infty w_n(t) \frac{z^n}{n!} = 
\frac{ \beta^2 }{ \beta^2 - z^2 } e^{tz},\qquad   t\in\mathbb{C};\  z\in\mathbb{C}:\ |z|<\beta. 
\end{equation}

\end{itemize}

\end{theorem}
\textbf{Proof.} 
In order to obtain the integral representation~(\ref{f3_78}) one can use Lagrange's method of variation of parameters.
We omit the details, since one
can check directly that the right-hand side of~(\ref{f3_78}) satisfies the required differential equation $D w_n = z^n$.
In fact, using the induction argument and the integration by parts one can verify that the improper integrals in~(\ref{f3_78})
exist and the right-hand side of~(\ref{f3_78}) is a monic polynomial of degree $n$.
Then the direct differentiation shows that $D w_n = z^n$. 
By Proposition~\ref{p2_1} the $n$-th degree polynomial solution of $D w = z^n$, is unique. Thus, 
the right-hand side of~(\ref{f3_78}) is indeed $w_n$.

Using the integration by parts two times we may write:
$$ \frac{2}{\beta} w_n(t) = 
\frac{2}{\beta} t^n +\frac{n}{\beta} e^{\beta t} \int_t^{+\infty} x^{n-1} e^{-\beta x} dx
- \frac{n}{\beta} e^{-\beta t} \int_{-\infty}^t x^{n-1} e^{\beta x} dx = $$
$$ = \frac{2}{\beta} t^n + \frac{n(n-1)}{\beta^2} 
\left(
e^{\beta t} \int_t^{+\infty} x^{n-2} e^{-\beta x} dx
+ e^{-\beta t} \int_{-\infty}^t x^{n-2} e^{\beta x} dx
\right) = $$
$$ = \frac{2}{\beta} t^n + \frac{2n(n-1)}{\beta^3} w_{n-2}(t),\qquad n\geq 2. $$ 
Therefore
\begin{equation}
\label{f3_81}
t^n = w_n(t) - \frac{n(n-1)}{\beta^2} w_{n-2}(t),\qquad n\geq 2.
\end{equation}
If we set $w_{-1}=w_{-2}=0$, then relation~(\ref{f3_81}) holds for all $n\in\mathbb{Z}_+$.
It remains to substitute expressions for $z^n$ and $z^{n+1}$ into 
$$ z^{n+1} = z z^n, $$
to obtain relation~(\ref{f3_79}). 

Denote by $Q(t,z)$ the function on the right in~(\ref{f3_80}).
It is an analytic function in the circle $|z|<\beta$.
Let us calculate its Taylor coefficients by the Leibniz rule:
$$ (Q(t,z))^{(k)}(0) = \left(
e^{tz} \frac{1}{1+\alpha z^2}
\right)^{(k)}(0) = $$
$$ = \sum_{j=0}^k \left(
\begin{array}{cc} k\\
j\end{array}\right)
\left(
\frac{1}{1+\alpha z^2}
\right)^{(j)}(0)
\left(
e^{tz} 
\right)^{(k-j)}(0) =
$$
$$ = \sum_{m=0}^{\left[ \frac{k}{2} \right]}
\left(
\begin{array}{cc} k\\
2m\end{array}\right)
(-\alpha)^m (2m)! t^{k-2m} = y_k(2,\alpha;t).
$$
The theorem is proved.
$\Box$

It is clear that relation~(\ref{f3_79}) can be written in the matrix form~(\ref{f1_43}). Thus, we conclude that
$w_n(z)$ are classical type Sobolev orthogonal polynomials on the unit circle.
Observe that
$w_2(z) = z^2-2\alpha$,
has roots
$z_{1,2} = \pm \sqrt{2|\alpha|} i$.
They are outside the unit circle, if $\alpha < -\frac{1}{2}$.
Thus, $w_n$ are not OPUC for $\alpha < -\frac{1}{2}$.

We remark that the generating function~$Q(t,z)$ from~(\ref{f3_80}) can be guessed from relation~(\ref{f3_81}), when dividing by $n!$ and
summing up formally. We came to this function in this way. Thus, this recurrence relation was very helpful.
In what follows we shall consider generating functions of a general form which include polynomials from Theorem~\ref{t3_1}
and many others. 

\begin{corollary}
\label{c3_2}
In conditions of Corollary~\ref{c3_1} with $r=2$, $\alpha<0$, the polynomials $\widehat y_n(2,\alpha;z)$ admit
the following integral representation:
$$ \widehat y_n(2,\alpha;z) = 
\frac{\beta}{2} e^{\beta t} 
\int_z^{+\infty} p_n(x) e^{-\beta x} dx  
+\frac{\beta}{2} e^{-\beta z}
\int_{-\infty}^{z} p_n(x) e^{\beta x} dx, $$ 
\begin{equation}
\label{f3_85}
z\in\mathbb{R},\ n\in\mathbb{Z}_+;\quad \beta:= \sqrt{ -\frac{1}{\alpha} }.  
\end{equation}

\end{corollary}
\textbf{Proof.}
Substitute into relation~(\ref{f3_73_2}) the integral representation~(\ref{f3_78}) and use~(\ref{f3_73_1}) to get relation~(\ref{f3_85}).
$\Box$

Let $\{ g_n(t) \}_{n=0}^\infty$ be a system of OPRL or OPUC, which has a generating function of the following form:
\begin{equation}
\label{f3_87}
G(t,w) = f(w) e^{ t u(w) } = \sum_{n=0}^\infty g_n(t) \frac{ w^n }{ n! },\qquad t\in\mathbb{C},\quad |w|< R_0,\quad (R_0>0),
\end{equation}
where $f,u$ are analytic functions in the circle $\{ |w|< R_0 \}$, $u(0)=0$.
Such generating functions for OPRL were studied by Meixner, see, e.g., \cite[p. 273]{cit_110_BE_3}.
As for the case of OPUC, we are not aware of any such a system, besides $\{ z^n \}_{n=0}^\infty$.
Consider the following function:
$$ F(t,w) = \frac{1}{ p(u(w)) } G(t,w) =  \frac{1}{ p(u(w)) } f(w) e^{ t u(w) },\qquad $$
\begin{equation}
\label{f3_89}
t\in\mathbb{C},\quad |w|< R_1 < R_0,\quad (R_1>0),
\end{equation}
where $p\in\mathbb{P}$: $p(0)\not = 0$. 
In the case $u(z)=z$, one should take $R_1\leq |z_0|$, where $z_0$ is a root of $p$ with the smallest modulus.  
This ensures that $F(t,w)$ is an analytic function with respect to each variable separately. Therefore $F$ is an analytic function
of two variables in any polydisk 
$C_{T_1,R_1} = \{ (t,w)\in\mathbb{C}^2:\ |t|< T_1,\ |w|< R_1 \}$, $T_1>0$.
In the general case, since $p(u(0))=p(0)\not=0$, there also exists a suitable $R_1$, which ensures that
$F$ is analytic in $C_{T_1,R_1}$.

Let us expand the function $F(t,w)$ in Taylor's series by $w$ with a fixed $t$:
\begin{equation}
\label{f3_91}
F(t,w) = \sum_{n=0}^\infty \varphi_n(t) \frac{ w^n }{ n! },\qquad (t,w)\in C_{T_1,R_1},
\end{equation}
where $\varphi_n(t)$ are some complex-valued functions.

\begin{lemma}
\label{l3_1}
The function $\varphi_n(t)$ is a complex polynomial of degree $n$, $\forall n\in\mathbb{Z}_+$.
\end{lemma}
\textbf{Proof.}
In fact,
$$ \varphi_n(t) = (F(t,w))^{(n)}_w|_{w=0} = 
\left.\left( \frac{1}{ p(u(w)) } G(t,w) \right) ^{(n)}_w \right|_{w=0} = $$
$$ = \sum_{j=0}^n \left(
\begin{array}{cc} n\\
j\end{array}\right)
\left(
\frac{1}{ p(u(w)) }
\right)^{(j)}(0)
\left.\left(
G(t,w)
\right)^{(n-j)}_w \right|_{w=0} = $$
\begin{equation}
\label{f3_91_1}
= \sum_{j=0}^n \left(
\begin{array}{cc} n\\
j\end{array}\right)
b_j
g_{n-j}(t),
\end{equation}
where
\begin{equation}
\label{f3_92}
b_j := \left(
\frac{1}{ p(u(w)) }
\right)^{(j)}(0). 
\end{equation}
Thus, $\varphi_n(t)$ is a polynomial. Since $b_0 = \frac{1}{ p(0) }\not=0$, then $\deg \varphi_n = n$.
$\Box$

Fix an arbitrary $k\in\mathbb{N}$. The function $F^{(k)}_t(t,w)$ is also analytic in $C_{T_1,R_1}$.
Consider its Taylor's expansion with respect to $w$, with a fixed $t$ :
\begin{equation}
\label{f3_93}
F^{(k)}_t(t,w) = \sum_{n=0}^\infty a_{n,k}(t) \frac{ w^n }{ n! },\qquad (t,w)\in C_{T_1,R_1},
\end{equation}
where $a_{n,k}(t)$ are some complex functions.

\begin{lemma}
\label{l3_2}
The functions $a_{n,k}(t)$ have the following form:
\begin{equation}
\label{f3_95}
a_{n,k}(t) = \varphi_n^{(k)}(t),\qquad n\in\mathbb{Z}_+,\ |t| < T_1;\ k\in\mathbb{N},
\end{equation}
i.e. we can differentiate the series~(\ref{f3_91}) with respect to $t$.
\end{lemma}
\textbf{Proof.}
We may write:
$$ a_{n,k}(t) = \left. \left( F^{(k)}_t \right)^{ (n) }_w \right|_{w=0} = 
\left. \left( F^{(n)}_w \right)^{ (k) }_t \right|_{w=0} = $$
$$ = \left.\left(
\varphi_n(t) + \sum_{l=n+1}^\infty \varphi_l(t) \frac{1}{ (l-n)! } w^{l-n}
\right) ^{ (k) }_t \right|_{w=0} = $$
$$ = \left.\left(
\varphi_n(t) + w L(t,w)
\right) ^{ (k) }_t \right|_{w=0}, $$
where 
$$ L(t,w) := \sum_{k=1}^\infty \varphi_{n+k}(t) \frac{ w^{k-1} }{k!},\qquad (t,w)\in C_{T_1,R_1}. $$
Since $\varphi_n(t) + w L(t,w)$, and $\varphi_n(t)$ (by Lemma~\ref{l3_1}) are analytic functions of two variables in $C_{T_1,R_1}$, 
then $w L(t,w)$ is also analytic in $C_{T_1,R_1}$. Therefore $L(t,w)$ is an analytic function of $t$ in $\{ |t|<T_1 \}$, for any fixed $w:\ |w|<R_1$.
By the definition of $L(t,w)$, it is represented by a power series with respect to $w$. Therefore $L(t,w)$ is an analytic
function of $w$ in $\{ |w|<R_1 \}$, for any fixed $t:\ |t|<T_1$.
Consequently, $L(t,w)$ is an analytic function of two variables in $C_{T_1,R_1}$.
Then
$$ a_{n,k}(t) =
\left.\left(
\varphi_n^{ (k) }(t) + w L^{ (k) }_t(t,w)
\right) \right|_{w=0} = \varphi_n^{ (k) }(t), $$
since $L^{ (k) }_t$ is an analytic function of $w$ in $\{ |w|<R_1 \}$, for any fixed $t$ in $\{ |t|<T_1 \}$.
$\Box$

Suppose that $\deg p\geq 1$, and
\begin{equation}
\label{f3_96}
p(z) = \sum_{k=0}^d c_k z^k,\qquad c_k\in\mathbb{C},\ c_d\not=0;\ c_0\not=0;\ d\in\mathbb{N}. 
\end{equation}
By~(\ref{f3_89}) we may write:
$$ \sum_{k=0}^d c_k F^{(k)}_t(t,w) = 
\sum_{k=0}^d c_k u^k(w) \frac{1}{ p(u(w)) } G(t,w)
= G(t,w),\qquad (t,w)\in C_{T_1,R_1}. $$
Comparing Taylor's coefficients of the series for both sides of the last relation (see~(\ref{f3_93}),(\ref{f3_95}) and (\ref{f3_87})) 
we conclude that
\begin{equation}
\label{f3_97}
\sum_{k=0}^d c_k \varphi_n^{(k)}(t) = g_n(t),\qquad n\in\mathbb{Z}_+.
\end{equation}
In particular, if $g_n(t)=t^n$, we obtain a generalization of an equation $D y(t)=t^n$, with the differential operator~(\ref{f3_30}),
which we have already studied.

\begin{theorem}
\label{t3_3}
Let $d\in\mathbb{N}$, and $p(z)$ be as in~(\ref{f3_96}).
Let $\{ g_n(t) \}_{n=0}^\infty$ be a system of OPRL or OPUC, having a generating function $G(t,w)$ from~(\ref{f3_87})
and $F(t,w)$ be given by~(\ref{f3_89}). Fix some positive $T_1, R_1$, such that $F(t,w)$ is analytic in the polydisk
$C_{T_1,R_1}$.
Polynomials 
\begin{equation}
\label{f3_100}
\varphi_n(z) = 
\sum_{j=0}^n \left(
\begin{array}{cc} n\\
j\end{array}\right)
b_j
g_{n-j}(t),\qquad n\in\mathbb{Z}_+,
\end{equation}
where $b_j = \left(
\frac{1}{ p(u(w)) }
\right)^{(j)}(0)$,
have the following properties:
\begin{itemize}
\item[(i)] Polynomials $\varphi_n$ are Sobolev orthogonal polynomials:
$$ \int \left( \varphi_n(t), \varphi_n'(t),..., \varphi_n^{(d)}(t) \right) \widetilde M \overline{
\left( \begin{array}{cccc} \varphi_m(t) \\
\varphi_m'(t) \\
\vdots \\
\varphi_m^{(d)}(t) \end{array} \right)
} 
d \mu_g = \tau_n \delta_{n,m},  $$
\begin{equation}
\label{f3_102}
\tau_n>0,\quad n,m\in\mathbb{Z}_+,
\end{equation}
where
\begin{equation}
\label{f3_104}
\widetilde M = (c_0,c_1,...,c_d)^T (\overline{c_0},\overline{c_1},...,\overline{c_d}).
\end{equation}
Here $d\mu_g$ is the measure of orthogonality of $g_n$.

\item[(ii)] Polynomials $\varphi_n$ have the generating function $F(t,w)$, and relation~(\ref{f3_91}) holds.

\item[(iii)] Polynomials $\varphi_n$ have the following integral representation:
\begin{equation}
\label{f3_106}
\varphi_n(t) = \frac{ n! }{2\pi i} \oint_{ |w|=R_2 }  
\frac{1}{ p(u(w)) } f(w) e^{ t u(w) } w^{-n-1} dw,\qquad n\in\mathbb{Z}_+,
\end{equation}
where $R_2$ is an arbitrary number, satisfying $0< R_2 <R_1$.

\end{itemize}

\end{theorem}
\textbf{Proof.} 
Representation~(\ref{f3_100}) for $\varphi_n$ was derived in the proof of Lemma~\ref{l3_1}.
Orthogonality relations for $\varphi_n$ follow from the orthogonality relations for $g_n$, by taking into account
relation~(\ref{f3_97}).
Finally, relation~(\ref{f3_106}) follows from the known representation of Taylor's coefficients of an analytic function.
$\Box$

There are two important cases of $g_n$, which lead to additional properties of $\varphi_n$, namely, to differential equations
and recurrence relations. 
Next two corollaries are devoted to these questions.

\begin{corollary}
\label{c3_3}
In conditions of Theorem~\ref{t3_3} suppose that $g_n(t) = t^n$, $n\in\mathbb{Z}_+$; $f(w)=1$, $u(w)=w$. 
Polynomials $\{ \varphi_n(t) \}_{n=0}^\infty$ satisfy the following recurrence relation:
$$ (n+1) \sum_{k=0}^d \varphi_{n+1-k}(t) \frac{c_k}{ (n+1-k)! } = $$
\begin{equation}
\label{f3_108}
= t \left(
\sum_{k=0}^d \varphi_{n-k}(t) \frac{c_k}{ (n-k)! }
\right),\qquad n\in\mathbb{Z}_+,
\end{equation}
where $\varphi_{r}:=0$, $r!:=1$, for $r\in\mathbb{Z}:\ r<0$.

Polynomials $\{ \varphi_n(t) \}_{n=0}^\infty$ obey the following differential equation:
\begin{equation}
\label{f3_110}
t \sum_{k=0}^d c_k \varphi_{n}^{(k+1)}(t)
= n \left(
\sum_{k=0}^d c_k \varphi_{n}^{ (k) }(t) 
\right),\qquad n\in\mathbb{Z}_+.
\end{equation}

\end{corollary}
\textbf{Proof.} In this case we have
$$ F(t,w) = \frac{1}{p(w)} e^{tw},\qquad (t,w)\in C_{T_1,R_1}. $$
Using power series we write:
$$ \sum_{k=0}^d \sum_{n=0}^\infty \varphi_n(t) c_k \frac{ w^{n+k} }{n!} = \sum_{n=0}^\infty t^n \frac{w^n}{n!},\qquad (t,w)\in C_{T_1,R_1}. $$
Changing the index $j=n+k$ on the left, and comparing Taylor's coefficients we obtain relation~(\ref{f3_108}).
In order to derive relation~(\ref{f3_108}) use relations~(\ref{f1_50}) and~(\ref{f3_97}).
$\Box$

\begin{corollary}
\label{c3_4}
In conditions of Theorem~\ref{t3_3} suppose that $g_n(t) = H_n(t)$, $n\in\mathbb{Z}_+$, are Hermite polynomials; $f(w)=e^{-w^2}$, $u(w)=2w$.
Polynomials $\{ \varphi_n(t) \}_{n=0}^\infty$ satisfy the following recurrence relation:
$$ (n+1) \sum_{k=0}^d \varphi_{n+1-k}(t) \frac{c_k 2^k}{ (n+1-k)! } +
2 \sum_{k=0}^d \varphi_{n-1-k}(t) \frac{c_k 2^k}{ (n+1-k)! }
= $$
\begin{equation}
\label{f3_112}
= 2t \left(
\sum_{k=0}^d \varphi_{n-k}(t) \frac{c_k 2^k}{ (n-k)! }
\right),\qquad n\in\mathbb{N},
\end{equation}
where $\varphi_{r}:=0$, $r!:=1$, for $r\in\mathbb{Z}:\ r<0$; and
\begin{equation}
\label{f3_112_1}
c_0 \varphi_1(t) + 2c_1 \varphi_{0}(t) = 2 c_0 t \varphi_0(t).
\end{equation}

Polynomials $\{ \varphi_n(t) \}_{n=0}^\infty$ obey the following differential equation:
\begin{equation}
\label{f3_114}
\sum_{k=0}^d c_k \varphi_{n}^{(k)}(t) - 2t \sum_{k=0}^d c_k \varphi_{n}^{(k+1)}(t)
= -2n \left(
\sum_{k=0}^d c_k \varphi_{n}^{ (k) }(t) 
\right),\qquad n\in\mathbb{Z}_+.
\end{equation}

\end{corollary}
\textbf{Proof.} In this case we have 
$G(t,w)=e^{2tw-w^2}$.
Rewriting
$$ p(2w) F(t,w) = G(t,w), $$
using power series of $w$, and comparing Taylor's coefficients, we get
$$ H_n(t) = n! \sum_{k=0}^d \varphi_{n-k}(t) c_k 2^k \frac{1}{(n-k)!},\qquad n\in\mathbb{Z}_+. $$
Using the three-term recurrence relation for the Hermite polynomials we obtain
relations~(\ref{f3_112}) and~(\ref{f3_112_1}).
By the differential equation:
$$ H_n''(t) - 2t H_n'(t) = -2n H_n(t), $$
and relation~(\ref{f3_97}) we obtain relation~(\ref{f3_114}).
$\Box$

Polynomials $\varphi_n$ from the last two corollaries fit into the scheme of Problem~1.
Notice that generating functions for $\varphi_n$ from these corollaries admit a study by the Darboux method,
see Theorem~8.4 from~\cite{cit_50000_Gabor_Szego}.
One should require that the zeros of $p$ are outside or on the unit circle (or use $p(w/c)$).

\begin{center}
{\large\bf 
On some Sobolev spaces with matrix weights and classical type Sobolev orthogonal polynomials.}
\end{center}
\begin{center}
{\bf S.M. Zagorodnyuk}
\end{center}

For every system $\{ p_n(z) \}_{n=0}^\infty$ of OPRL or OPUC, we construct Sobolev orthogonal polynomials $y_n(z)$, 
with explicit integral representations
involving $p_n$. 
Two concrete families of Sobolev orthogonal polynomials (depending on an arbitrary number of
complex parameters) which are generalized eigenvalues of a difference operator (in $n$)
and generalized eigenvalues of a differential operator (in $n$) are given.
Applications of a general connection between Sobolev orthogonal polynomials and 
orthogonal systems of functions in the direct sum of scalar $L^2_\mu$ spaces
are discussed.

\vspace{1.5cm}

V. N. Karazin Kharkiv National University \newline\indent
School of Mathematics and Computer Sciences \newline\indent
Department of Higher Mathematics and Informatics \newline\indent
Svobody Square 4, 61022, Kharkiv, Ukraine

Sergey.M.Zagorodnyuk@gmail.com; Sergey.M.Zagorodnyuk@univer.kharkov.ua

}

\begin{thebibliography}{99}

\bibitem{cit_5_Azad}
Azad H.; Laradji A.; Mustafa M. T. Polynomial solutions of differential equations. Adv. Difference Equ. 2011:58 (2011), 12 pp.

\bibitem{cit_100_Duran}
Duran, Antonio J. A generalization of Favard's theorem for polynomials satisfying a recurrence relation. 
J. Approx. Theory 74 (1993), no. 1, 83--109.

\bibitem{cit_110_BE_3}
Erd\'elyi, Arthur; Magnus, Wilhelm; Oberhettinger, Fritz; Tricomi, Francesco G. 
Higher transcendental functions. Vol. III. Based, in part, on notes left by Harry Bateman. 
McGraw-Hill Book Company, Inc., New York-Toronto-London, 1955. {\rm xvii}+292 pp.


\bibitem{cit_1500___Kim__2014}
Kim, H. K.; Kwon, K. H.; Littlejohn, L. L.; Yoon, G. J. 
Diagonalizability and symmetrizability of Sobolev-type bilinear forms: a combinatorial approach. 
Linear Algebra Appl. 460 (2014), 111--124.  

\bibitem{cit_1700__Klotz_1992}
Klotz, L. A matrix generalization of a theorem of Szeg\H{o}. Anal. Math. 18 (1992), no. 1, 63--72.

\bibitem{cit_5105_Koekoek_book}
Koekoek R., Lesky P. A., Swarttouw R. F. Hypergeometric orthogonal polynomials and their 
$q$-analogues. With a foreword by Tom H. Koornwinder. Springer Monographs in Mathematics. Springer-Verlag, Berlin, (2010).


\bibitem{cit_2000__Kwon___2009}
Kwon, K. H.; Littlejohn, Lance L.; Yoon, G. J. Ghost matrices and a characterization of 
symmetric Sobolev bilinear forms. Linear Algebra Appl. 431 (2009), no. 1-2, 104--119. 


\bibitem{cit_5150_Survey_2006}
Marcell\'an, Francisco; Moreno Balc\'azar, Juan Jos\'e. Asymptotics and zeros of Sobolev orthogonal polynomials 
on unbounded supports. Acta Appl. Math. 94 (2006), no. 2, 163--192.  


\bibitem{cit_5150_M_X_Survey_2015}
Marcell\'an, Francisco; Xu, Yuan. On Sobolev orthogonal polynomials. Expo. Math. 33 (2015), no. 3, 308--352.

\bibitem{cit_5150_Survey_1993}
Marcell\'an, F.; Alfaro, M.; Rezola, M. L. Orthogonal polynomials on Sobolev spaces: old and new directions. 
Proceedings of the Seventh Spanish Symposium on Orthogonal Polynomials and Applications (VII SPOA) (Granada, 1991). 
J. Comput. Appl. Math. 48 (1993), no. 1-2, 113--131.  

\bibitem{cit_7000_Markus}
Markus A. S. Introduction to the spectral theory of polynomial operator pencils. 
With an appendix by M. V. Keldysh. Translations of Mathematical Monographs, 71. American Mathematical Society, 
Providence, RI, (1988).

\bibitem{cit_8000_Rodman}
Rodman L. An introduction to operator polynomials. 
Operator Theory: Advances and Applications, 38. Birkh\"auser Verlag, Basel, (1989).

\bibitem{cit_8100_Marden}
Marden, Morris. Geometry of polynomials. Second edition. Mathematical Surveys, No. 3 American Mathematical Society, Providence, R.I. 
1966 {\rm xiii}+243 pp.

\bibitem{cit_5150__Survey_1998}
Mart\'inez-Finkelshtein, Andrei. Asymptotic properties of Sobolev orthogonal polynomials. 
Proceedings of the VIIIth Symposium on Orthogonal Polynomials and Their Applications (Seville, 1997). 
J. Comput. Appl. Math. 99 (1998), no. 1-2, 491--510. MR1662717 Add to clipboard 

\bibitem{cit_5150__Survey_2001}
Mart\'inez-Finkelshtein, A. Analytic aspects of Sobolev orthogonal polynomials revisited. Numerical analysis 2000, Vol. V, 
Quadrature and orthogonal polynomials. J. Comput. Appl. Math. 127 (2001), no. 1-2, 255--266. 

\bibitem{cit_5150__Survey_1996}
Meijer, H. G. A short history of orthogonal polynomials in a Sobolev space. 
I. The non-discrete case. 31st Dutch Mathematical Conference (Groningen, 1995). 
Nieuw Arch. Wisk. (4) 14 (1996), no. 1, 93--112.

\bibitem{cit_7000__Rosenberg_1964}
Rosenberg, Milton. The square-integrability of matrix-valued functions with respect to a non-negative Hermitian measure. 
Duke Math. J. 31 (1964), 291--298.


\bibitem{cit_50000_Gabor_Szego}
Szeg\"o, G\'abor. Orthogonal polynomials. Fourth edition. 
American Mathematical Society, Colloquium Publications, Vol. XXIII. American Mathematical Society, Providence, R.I., 1975. xiii+432 pp.


\bibitem{cit_8000__Weron_1974}
Weron, A. On characterizations of interpolable and minimal stationary processes. Studia Math. 49 (1973/74), 165--183.

\bibitem{cit_9000__Zagorodnyuk_JAT_2020}
Zagorodnyuk, S.M. On some classical type Sobolev orthogonal polynomials. J. Approx. Theory 250 (2020), 105337.

\bibitem{cit_9000__Zagorodnyuk_CMA_2020}
Zagorodnyuk, S.M. On a family of hypergeometric Sobolev orthogonal polynomials on the unit circle. Constr. Math. Anal.
3(2) (2020), 84--75.

\bibitem{cit_9200__Xu_2018}
Xu, Yuan. Approximation by polynomials in Sobolev spaces with Jacobi weight. J. Fourier Anal. Appl. 24 (2018), no. 6, 1438--1459.

\end{thebibliography}
\end{document}